# $L_1$-norm minimization for quaternion signals

Jiasong Wu, *Member, IEEE*, Xu Zhang, Xiaoqing Wang, Lotfi Senhadji, *Senior Member, IEEE*, and Huazhong Shu, *Senior Member, IEEE*

*Abstract*—The $l_1$-norm minimization problem plays an important role in the compressed sensing (CS) theory. We present in this letter an algorithm for solving the problem of $l_1$-norm minimization for quaternion signals by converting it to second-order cone programming. An application example of the proposed algorithm is also given for practical guidelines of perfect recovery of quaternion signals. The proposed algorithm may find its potential application when CS theory meets the quaternion signal processing.

## I. INTRODUCTION

THE problem of $l_1$-norm minimization plays an important role in the recently developed compressed sensing (CS) theory [1]-[3], which is a new approach for data acquisition that finds wide applications in the field of signal and image processing. CS has been conventionally used to real input data. Recently, special attention has also been paid to the complex input data for addressing problems such as blind source separation [4], radar signal processing [5], and terahertz imaging [6].

On the other hand, the theory and application of quaternion (or hypercomplex) algebra, invented by Hamilton [7], has received many attentions in recent years (for example, [8]-[22]). Ell and Sangwine [8] applied the quaternion Fourier transform to color image processing. Sangwine and Bihan [9] then developed the Quaternion Toolbox for Matlab (QTFM). Jiang and Wei [16] proposed an algorithm for solving the quaternion least squares problem:

$$\min \quad \|\mathbf{Ax}-\mathbf{y}\|_2, \quad \text{s.t.} \quad \mathbf{Bx}=\mathbf{z}, \qquad (1)$$

where $\mathbf{A}\in\mathbb{Q}^{n\times m}$, $\mathbf{B}\in\mathbb{Q}^{u\times m}$, $\mathbf{x}\in\mathbb{Q}^{m\times 1}$, $\mathbf{y}\in\mathbb{Q}^{n\times 1}$, $\mathbf{z}\in\mathbb{Q}^{u\times 1}$, $\mathbb{Q}$ denotes the quaternion number field, and $\|\cdot\|_2$ denotes the $l_2$-norm of quaternion vector. "s.t." means subject to. The algorithm reported in [16] was further extended by Jiang *et al.* [17] and Wang *et al.* [18] to two-dimensional where the inputs are quaternion matrices. Note that the algorithms reported in [16]-[18] solve the overdetermined system, that is, $n\geq m$ for measurement matrix $\mathbf{A}$.

To the authors' knowledge, the $l_1$-norm minimization problem for quaternion signal has not yet been investigated. In [4], Winter *et al.* converted the $l_1$-norm minimization of complex signal to second-order cone programming (SOCP), which was then solved by SeDuMi software [23]. In this paper, we extend the algorithm in [4] to the case where the input data are quaternion signals and we derive a new algorithm for solving the following optimization problem:

$$\min \quad \|\mathbf{x}\|_1, \quad \text{s.t.} \quad \mathbf{y}=\mathbf{Ax}, \qquad (2)$$

where $\mathbf{A}$, $\mathbf{x}$, and $\mathbf{y}$ are defined in (1). $\|\cdot\|_1$ denotes the $l_1$-norm of quaternion vector. Here, we assume that $n\leq m$, that is, we deal with the underdetermined system for measurement matrix $\mathbf{A}$, which is different from that of (1). The proposed algorithm may find its potential applications when CS theory meets the quaternion signal processing (for example, quaternion blind source separation [19]-[21]).

The rest of this letter is organized as follows. Section II recalls some preliminaries about the quaternions. Section III describes the proposed algorithm for $l_1$-norm minimization problem of quaternion signals. An application example of the proposed algorithm is provided in Section IV for practical guidelines. Section V concludes the paper and also points out the direction of further research.

## II. DEFINITIONS

A quaternion is a hypercomplex number which consists of one real part and three imaginary parts as follows:

$$q = R(q)+I(q)i+J(q)j+K(q)k, \qquad (3)$$

where $R(q)$, $I(q)$, $J(q)$, $K(q)\in\mathbb{R}$, $\mathbb{R}$ denotes the real number field, and *i*, *j*, and *k* are three imaginary units obeying the following rules:

$$i^2 = j^2 = k^2 = -1, \qquad (4)$$

This work was supported by the National Basic Research Program of China under Grant 2011CB707904, by the National Natural Science Foundation of China under Grants 60873048 and 60911130370.

J. Wu is with the Laboratory of Image Science and Technology, School of Computer Science and Engineering, Southeast University, Nanjing 210096, China, the Centre de Recherche en Information Biomédicale Sino-Français (CRIBs), Nanjing 210096, China, INSERM U 1099, 35000 Rennes, France, and the Laboratoire Traitement du Signal et de l'Image (LTSI), Université de Rennes 1, 35000 Rennes, France (e-mail: jswu@seu.edu.cn).

X. Zhang, X. Wang and H. Shu are with the Laboratory of Image Science and Technology, School of Computer Science and Engineering, Southeast University, Nanjing 210096, China, and also with the Centre de Recherche en Information Biomédicale Sino-Français (CRIBs), Nanjing 210096, China (e-mail: zzal119911@gmail.com; xiaoqingwang2010@gmail.com; shu.list@seu.edu.cn).

L. Senhadji is with INSERM U1099, 35000 Rennes, France, with the Laboratoire Traitement du Signal et de l'Image (LTSI), Université de Rennes 1, 35000 Rennes, France, and with the Centre de Recherche en Information Biomédicale Sino–Français (CRIBs), 35000 Rennes, France (e-mail: lotfi.senhadji@univ-rennes1.fr).

$$ij = -ji = k, \ jk = -kj = i, \ ki = -ik = j. \quad (5)$$

The conjugate and modulus of a quaternion are, respectively, defined by

$$q^* = R(q) - I(q)i - J(q)j - K(q)k, \quad (6)$$

$$|q| = \sqrt{R^2(q) + I^2(q) + J^2(q) + K^2(q)}. \quad (7)$$

We also define

$$\mathbf{A} = [\mathbf{a}_1, \cdots, \mathbf{a}_m], \quad (8)$$

$$\mathbf{x} = [x_1, \cdots, x_m]^T = R(\mathbf{x}) + I(\mathbf{x})i + J(\mathbf{x})j + K(\mathbf{x})k, \quad (9)$$

$$\mathbf{y} = [y_1, \cdots, y_n]^T = R(\mathbf{y}) + I(\mathbf{y})i + J(\mathbf{y})j + K(\mathbf{y})k. \quad (10)$$

The $l_p$-norm of a quaternion vector $\mathbf{x}$ is given by

$$\|\mathbf{x}\|_p = \left(\sum_{r=1}^m |x_r|^p\right)^{1/p}, \ p = 1, 2. \quad (11)$$

## III. ALGORITHM

In this section, we derive an algorithm for $l_1$-norm minimization problem with quaternion signals by using SOCP. Equation (2) is equivalent to

$$\min \ t \in \mathbb{R}^+, \quad \text{s.t.} \ \mathbf{y} = \mathbf{A}\mathbf{x}, \quad \|\mathbf{x}\|_1 \leq t. \quad (12)$$

By decomposing $t = \sum_{r=1}^m t_r$, $t_r \in \mathbb{R}^+$, $\mathbb{R}^+$ denotes positive real number field, the second constraint $\|\mathbf{x}\|_1 \leq t$ can be expressed as

$$\begin{aligned}
\|\mathbf{x}\|_1 &= \sum_{r=1}^m |x_r| \\
&= \sum_{r=1}^m \sqrt{R^2(x_r) + I^2(x_r) + J^2(x_r) + K^2(x_r)} \\
&= \sum_{r=1}^m \left\| \begin{bmatrix} R(x_r) \\ I(x_r) \\ J(x_r) \\ K(x_r) \end{bmatrix} \right\|_2 \leq \mathbf{1}^T \mathbf{t} = \mathbf{1}^T [t_1, t_2, \cdots, t_m]^T = t.
\end{aligned} \quad (13)$$

Then, (12) becomes

$$\min_{\mathbf{t}} \ \mathbf{1}^T \mathbf{t} \in \mathbb{R},$$

$$\text{s.t. } \mathbf{y} = \mathbf{A}\mathbf{x}, \quad \left\| \begin{bmatrix} R(x_r) \\ I(x_r) \\ J(x_r) \\ K(x_r) \end{bmatrix} \right\|_2 \leq t_r, \ \forall r \in \{1, 2, \cdots, m\}. \quad (14)$$

After some manipulation, equation (14) can be written as

$$\min_{\hat{\mathbf{x}}} \ \hat{\mathbf{c}}^T \hat{\mathbf{x}} \in \mathbb{R},$$

$$\text{s.t. } \hat{\mathbf{y}} = \hat{\mathbf{A}} \hat{\mathbf{x}}, \quad \left\| \begin{bmatrix} R(x_r) \\ I(x_r) \\ J(x_r) \\ K(x_r) \end{bmatrix} \right\|_2 \leq t_r, \ \forall r \in \{1, 2, \cdots, m\}, \quad (15)$$

where

$$\hat{\mathbf{x}} = [t_1, R(x_1), I(x_1), J(x_1), K(x_1), \cdots, \\ t_m, R(x_m), I(x_m), J(x_m), K(x_m)]^T \in \mathbb{R}^{5m}, \quad (16)$$

$$\hat{\mathbf{c}} = [1, 0, 0, 0, 0, \cdots, 1, 0, 0, 0, 0]^T \in \mathbb{R}^{5m}, \quad (17)$$

$$\hat{\mathbf{y}} = [R(\mathbf{y}), I(\mathbf{y}), J(\mathbf{y}), K(\mathbf{y})]^T \in \mathbb{R}^{4n}, \quad (18)$$

$$\hat{\mathbf{A}} = \begin{bmatrix}
0 & R(\mathbf{a}_1) & -I(\mathbf{a}_1) & -J(\mathbf{a}_1) & -K(\mathbf{a}_1) & \cdots \\
0 & I(\mathbf{a}_1) & R(\mathbf{a}_1) & -K(\mathbf{a}_1) & J(\mathbf{a}_1) & \cdots \\
0 & J(\mathbf{a}_1) & K(\mathbf{a}_1) & R(\mathbf{a}_1) & -I(\mathbf{a}_1) & \cdots \\
0 & K(\mathbf{a}_1) & -J(\mathbf{a}_1) & I(\mathbf{a}_1) & R(\mathbf{a}_1) & \cdots \\
& 0 & R(\mathbf{a}_m) & -I(\mathbf{a}_m) & -J(\mathbf{a}_m) & -K(\mathbf{a}_m) \\
& 0 & I(\mathbf{a}_m) & R(\mathbf{a}_m) & -K(\mathbf{a}_m) & J(\mathbf{a}_m) \\
& 0 & J(\mathbf{a}_m) & K(\mathbf{a}_m) & R(\mathbf{a}_m) & -I(\mathbf{a}_m) \\
& 0 & K(\mathbf{a}_m) & -J(\mathbf{a}_m) & I(\mathbf{a}_m) & R(\mathbf{a}_m)
\end{bmatrix} \in \mathbb{R}^{4n \times 5m}, \quad (19)$$

Equation (15) is the standard form of SOCP problem and can be solved by using several mature toolboxes, such as SeDuMi [23]. Then, we can easily obtain the recovered quaternion signal **xr** from $\hat{\mathbf{x}}$.

## IV. AN APPLICATION EXAMPLE

In this section, just like [1], we present numerical experiments that suggest empirical bounds on sparsity $s$ (time domain support of the input signal) relative to $n$ (the number of measurements) for a quaternion signal $\mathbf{x}$ to be perfectly recovered. The results can be seen as a set of practical guidelines for situations where one can expect perfect recovery from random Gaussian quaternion measurement matrix information using SOCP.

The following program is carried out on a standard desktop computer PC machine, which sets up Microsoft Windows 7 ultimate version system with Intel Core 2 Duo Processor E8400 (3 GHz) and 2 GB RAM. The algorithm has been implemented by Matlab R2011b.

Our experiments are as follows:
1) Random Gaussian measurement matrix $\mathbf{A} \in \mathbb{Q}^{n \times m} (n \leq m)$: choose constants $m = 256$ (the length of input signal), $n$ (the number of measurements), and then generate $\mathbf{A}$ with random entries sampled from independent and identically distributed (i.i.d.) Gaussian process with zero mean and variance equal to 1 (in quaternion $l_2$-norm sense).
2) Sparse quaternion input signal $\mathbf{x} \in \mathbb{Q}^{m \times 1}$: select a support set $\mathbf{T}$ of size $|\mathbf{T}| = s$ (sparsity) uniformly at random, and sample a vector $\mathbf{x}$ on $\mathbf{T}$ with i.i.d. Gaussian entries. For more discussion about quaternion Gaussian vector, we refer the readers to [11] and [12].
3) Quaternion output signal $\mathbf{y} \in \mathbb{Q}^{n \times 1}$: multiplied $\mathbf{A}$ by input $\mathbf{x}$ to obtain $\mathbf{y}$.
4) Construct the vectors $\hat{\mathbf{x}}$, $\hat{\mathbf{c}}$, $\hat{\mathbf{y}}$, and matrix $\hat{\mathbf{A}}$ (see from (16) to (19)).
5) Call the SeDuMi toolbox [23] to solve SOCP problem in (15) and compute the error, which is the $l_2$-norm of the deference between the recovered signal **xr** and the input signal $\mathbf{x}$, that is, $\|\mathbf{xr} - \mathbf{x}\|_2$.

6) We perform 100 times experiments for each pair of $s$ and $n$ and then save these errors and also count the number of perfect recovered experiment. The decision for perfectly recovery is $\|\mathbf{xr} - \mathbf{x}\|_2 \leq 10^{-8}$.

Fig. 1(a) shows the number of measurements $n$ with respect to sparsity $s$ for $m = 256$. Fig. 1(b) gives the cross section of the image in (a) at $n = 8, 16, 32, 64, 128$. It can be seen from this figure that for $n \geq 32$, the recovery rate is about 80% when $s \leq n/3$ and practically 100% when $s \leq n/4$. The result is better than that of the input data being real signal [1], whose recovery rate is about 80% when $s \leq n/5$ and practically 100% when $s \leq n/8$. We think the reason is that every quaternion signal value consist of 4 real values, so, the sparsity of quaternion signal is similar to the block sparsity (the length of block is 4) in real signal, which can be seen as a special sparsity form.

To further illustrate the recovery experiment results, Fig. 2 depicts the original quaternion signal $\mathbf{x}$ and its corresponding recovery signal $\mathbf{xr}$. In Fig. 2(a), $m = 256$, $n = 128$, and $s = 60$ and the quaternion signal is perfectly recovered. In Fig. 2(b), $m = 256$, $n = 128$, and $s = 90$ and in this case the recovery process failed.

## V. CONCLUSION

In this letter, we derived an algorithm for solving the $l_1$-norm minimization problem of quaternion signals, which is converted to SOCP and then solved by SeDuMi software. A numerical example was provided to illustrate the feasibility of the algorithm. The results can be viewed as a set of practical guidelines for situations where one can expect perfect recovery from random Gaussian quaternion measurement matrix information using SOCP. Further research includes: 1) study the robustness of the recovery procedure from the incomplete and inaccurate measurements of quaternion compressed sensing (QCS) as that of in [3]; 2) extend the quaternion signal vector $l_1$-norm minimization algorithm to that of quaternion matrix nuclear norm minimization, that is, quaternion matrix completion (QMC) [24]; 3) study the problem of quaternion robust principal component analysis (QRPCA) [22], [25]; 4) study the Bayesian framework of QCS, QMC, and QRPCA [26]-[28].

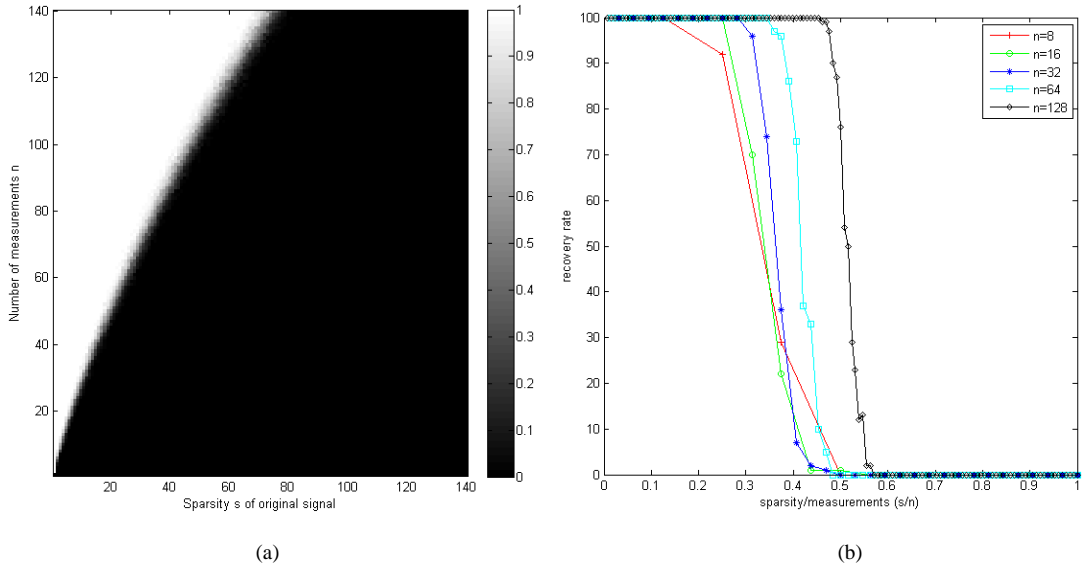

Fig. 1. Recovery experiment for *m*=256. (a) The image intensity standards for the percentage of the time solving (15) recovered the quaternion signal **x** as a function of the number of measurements *n* (vertical axis) and sparsity *s* (horizontal axis); the signal is recovered approximately 100% in white regions and never recovered in black regions. (b) Cross section of the image in (a) at *n* = 8, 16, 32, 64, 128.

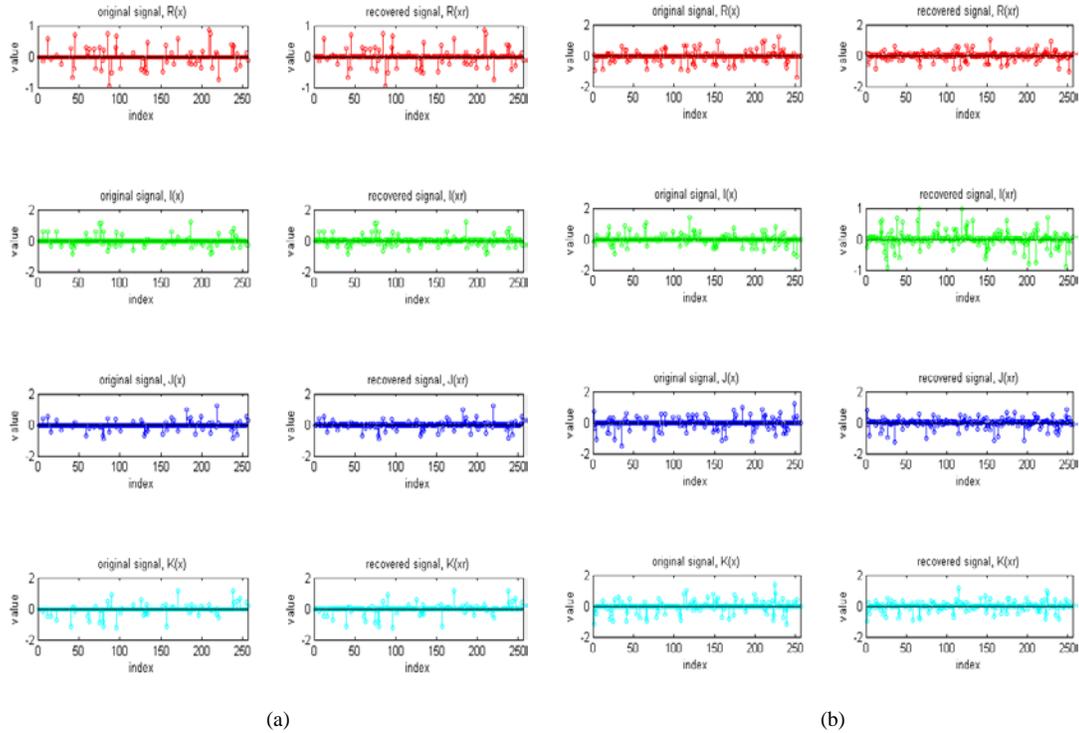

Fig. 2. Recovery experiment for $l_1$-norm minimization for quaternion signals (*m*=256 and *n*=128). (a) Left: original quaternion signal **x** with *s*=60; Right: the perfectly recovered signal **xr**. (b) Left: original quaternion signal **x** with *s*=90; Right: the failed recovered signal **xr**.